\documentclass{article}

\usepackage[english]{babel}

\usepackage[letterpaper,top=2cm,bottom=2cm,left=3cm,right=3cm,marginparwidth=1.75cm]{geometry}

\usepackage{amssymb}
\usepackage{latexsym}
\usepackage{dsfont}
\usepackage{amsmath}
\usepackage{graphicx}
\usepackage[colorlinks=true, allcolors=blue]{hyperref}

\allowdisplaybreaks[4]               

\newtheorem{thm}{Theorem}[section]

\newtheorem{lem}[thm]{Lemma}

\numberwithin{equation}{section}

\newcommand{\Cov}{{\mathbb{C}}{\rm ov}}
\newcommand{\Var}{{\mathbb{V}}{\rm ar}}

\newcommand{\Enum}{\mathbb{E}}

\newcommand{\Pnum}{\mathbb{P}}

\newcommand{\qed}{\hfill $\Box$}

\title{On the First and the Second Borel-Cantelli Lemmas}
\author{Jian-Sheng Xie\footnote{Corresponding author. Email: jsxie@fudan.edu.cn.} \, and Qihang Wang\\
{\footnotesize School of Mathematical Sciences, Fudan
University, Shanghai 200433,   China}} 
\date{}

\begin{document}
\maketitle

\begin{abstract}
Let $\{A_n\}_{n=1}^\infty$ be a sequence of events and let $\displaystyle S:=\sum_{n=1}^\infty 1_{A_n}$. We present in this note equivalent characterizations for
the statements $\Pnum (S<\infty)=1$ and $\Pnum (S=\infty)=1$ respectively. These characterizations are of Borel-Cantelli lemma type and of Kochen-Stone lemma type
respectively, which could be regarded as the most general version of the first and the second Borel-Cantelli Lemmas.

\noindent\textbf{MSC(2020):} 60A05; 60A10; 60F15;
\end{abstract}

\section{Introduction}\label{sec: C1}
The first and the second Borel-Cantelli lemmas are foundamental tools in discussing almost sure convergence. In this note we will present a version of their
generalizations, namely, Theorems \ref{thm: main-A} and \ref{thm: main-B}. 

The following is the so called Borel-Cantelli lemma (also called the first Borel-Cantelli lemma).
\begin{lem} (\textbf{Borel-Cantelli Lemma})
If a sequence $\{A_n\}_{n=1}^\infty$ of events satisfies
\begin{equation}\label{eq: prob-sum-finite}
\sum_{n=1}^\infty \Pnum (A_n)<\infty,
\end{equation}
then one has
\begin{equation}\label{eq: a.s.-finite}
\sum_{n=1}^\infty 1_{A_n}<\infty \quad \hbox{almost surely.}
\end{equation}
\end{lem}
In this note we prove the following theorem characterizing \eqref{eq: a.s.-finite}, which is of Borel-Cantelli lemma type.
\begin{thm}\label{thm: main-A}
Eq. \eqref{eq: a.s.-finite} holds true if and only if there exists an increasing subsequence $\{n_k\}_{k=1}^\infty$ of natural numbers such that
\begin{equation}\label{eq: prob-sum-B-finite}
\sum_{k=1}^\infty \Pnum (B_k)<\infty,
\end{equation}
where $\displaystyle B_k:=\bigcup_{j \in (n_{k-1}, n_k]} A_j$ with the convention $n_0:=0$.
\end{thm}

There are several versions of the second Borel-Cantelli lemma; The following is a synthesis.
\begin{lem} (\textbf{2nd Borel-Cantelli Lemma})
Given a sequence of events $\{A_n\}_{n=1}^\infty$ with
\begin{equation}\label{eq: prob-sum-infinite}
\sum_{n=1}^\infty \Pnum (A_n)=\infty.
\end{equation}
If they satisfy one of the following conditions
\begin{enumerate}
  \item[{\rm (1)}] The events $\{A_n\}_{n=1}^\infty$ are independent to each other;
  \item[{\rm (2)}] The events $\{A_n\}_{n=1}^\infty$ are pairwise-independent;
  \item[{\rm (3)}] There exists an $\ell^1$ sequence $\{\rho_{_k}\}_{k=1}^\infty$ of positive real numbers such that
$$
\Pnum (A_i \cap A_j) -\Pnum (A_i)\Pnum (A_j) \leq \rho_{_{|i-j|}} \cdot [\Pnum (A_i)+\Pnum (A_j)], \quad \forall i \neq j;
$$
  \item[{\rm (4)}] There exists an infinite matrix $M=(M_{_{i,j}})_{_{1 \leq i, j <\infty}}$ which induces a bouned continuous linear transformation
$$
\ell^2 \to \ell^2, x=(x_1, x_2, \cdots) \mapsto x M
$$
such that
$$
\Pnum (A_i \cap A_j) -\Pnum (A_i)\Pnum (A_j) \leq M_{_{i,j}} \cdot \sqrt{\Pnum (A_i)\Pnum (A_j)}, \quad \forall i \neq j;
$$
\end{enumerate}
then one has
\begin{equation}\label{eq: a.s.-infinite}
\sum_{n=1}^\infty 1_{A_n}=\infty \quad \hbox{almost surely.}
\end{equation}
\end{lem}
The above lemma with condition {\rm (1)} is the classic one which appears in most probability textbooks. The other conditions' versions are scattered in literatures:
\cite[Theorem 2.3.9]{Durrett} is for condition {\rm (2)}; \cite[Theorem 3.2]{Simon} is for condition {\rm (4)} (traced back to Erd\"{o}s and Rennyi \cite{ER}); condition 
{\rm (3)} is recorded in \cite{DMR2022} as uniform mixing condition. And recently we have established the following result which can be used in proving the above lemma 
in a unified manner.
\begin{thm}\cite{XZ}
Let $\{S_n\}_{n=1}^\infty$ be an increasing sequence of positive random variables. Suppose $\Enum S_n \to \infty$ and  there exists positive constants $C, \delta>0$ such 
that
\begin{equation}
\Var (S_n) \leq C(\Enum S_n)^{2 -\delta}, \quad \forall n \geq 1,
\end{equation}
or even more weakly
\begin{equation}
\Var (S_n) \leq \frac{C(\Enum S_n)^2}{(\ln \Enum S_n)^{1+\delta}}, \quad \forall n \geq 1,
\end{equation}
then
\begin{equation}
\lim_{n \to \infty} \frac{S_n}{\Enum S_n}=1  \quad \hbox{almost surely.}
\end{equation}
\end{thm}

But a more general version of the second Borel-Cantelli lemma in literature is the following one due to Kochen and Stone \cite{KS1964} \cite[Exercise 2.3.10]{Durrett}.
\begin{lem} (\textbf{Kochen-Stone Lemma})
Given a sequence of events $\{A_n\}_{n=1}^\infty$ with
\begin{equation}\label{eq: prob-sum-infinite}
\sum_{n=1}^\infty \Pnum (A_n)=\infty.
\end{equation}
Let $S_n:=1_{A_1}+\cdots+1_{A_n}$. If they satisfy
\begin{equation}\label{eq: limit4Var-Expect}
\varliminf_{n \to \infty} \frac{\Var (S_n)}{(\Enum S_n)^2}=0,
\end{equation}
then \eqref{eq: a.s.-infinite} holds true.
\end{lem}
In this note we prove the following theorem characterizing \eqref{eq: a.s.-infinite}, still of Kochen-Stone lemma type.
\begin{thm}\label{thm: main-B} 
Eq. \eqref{eq: a.s.-infinite} holds true if and only if there exists an increasing subsequence $\{n_k\}_{k=1}^\infty$ of natural numbers such that
\begin{equation}\label{eq: prob-sum-B-infinite}
\sum_{k=1}^\infty \Pnum (B_k)=\infty
\end{equation}
and \eqref{eq: limit4Var-Expect} hold true with $\displaystyle B_k:=\bigcup_{j \in (n_{k-1}, n_k]} A_j$ and $S_m:=1_{B_1}+\cdots+1_{B_m}$.
\end{thm}

More than a century passed ever since the establishments of the first and the second Borel-Cantelli lemmas (Borel, 1909; Cantelli, 1917; Cf. \cite{EH2022}). It's somehow weird 
that, as to our knowledge, no account of our theorems \ref{thm: main-A} and \ref{thm: main-B} can be found in the existing literatures, neither in all classical textbooks of 
probability theory, nor in other published literatures such as the monographs \cite{Chandra2012} \cite{DMR2022} or research papers (e.g., \cite{CE1952} \cite{Lamperti1963} 
\cite{DF1965} \cite{Freedman1973} \cite{Davie2001} \cite{Hoover2002} \cite{Petrov2002} \cite{Luzia2014} \cite{Frolov2021} \cite{EH2022} and references therein) on Borel-Cantelli 
lemmas' topic. Though currently the two theorems are of not much use in applications, they are theoretically significant in proving that the essential arguments in discussing 
almost sure convergence can always be regarded as Borel-Cantelli lemma and Kochen-Stone lemma. Noting also that, there is a trend of developping dynamical Borel–Cantelli lemmas in 
dynamical systems (relating especially to  number theory) ever since the publication of \cite{Philipp1967} (Cf. \cite{CK2001} \cite{Kifer2021} \cite{DMR2022} \cite{HLSW2022}), our 
theorems may have possible appliactions there.

In a word, Theorems \ref{thm: main-A} and \ref{thm: main-B} can be viewed respectively as the most general version of the first and the second Borel-Cantelli lemmas. And we will prove them in the next section.
\section{Proof of Theorems \ref{thm: main-A} and \ref{thm: main-B}}\label{sec: C2}
Given an increasing subsquence $\{n_k\}_{k=1}^\infty$ of natural numbers and a sequence $\{A_n\}_{n =1}^\infty$ of events, define $\displaystyle B_k:=\bigcup_{j \in (n_{k-1}, n_k]} A_j$ with
convention $n_0:=0$. It is easy to see that
\begin{eqnarray}
\label{eq: equvi-finite} \sum_{n=1}^\infty 1_{A_n} <\infty &\Leftrightarrow& \sum_{k=1}^\infty 1_{B_k} <\infty,\\
\label{eq: equvi-infinite} \sum_{n=1}^\infty 1_{A_n} =\infty &\Leftrightarrow& \sum_{k=1}^\infty 1_{B_k} =\infty.
\end{eqnarray}

\noindent \textbf{Proof of Theorem \ref{thm: main-A}.\;} Suppose \eqref{eq: prob-sum-B-finite} holds true for some increasing subsquence $\{n_k\}_{k=1}^\infty$ of natural numbers. Borel-Cantelli lemma
tells us $\displaystyle \sum_{k=1}^\infty 1_{B_k} <\infty$ almost surely. Then \eqref{eq: equvi-finite} tells us the validity of \eqref{eq: a.s.-finite}.

Now assume the validity of \eqref{eq: a.s.-finite}. It is just $\displaystyle \Pnum (\varlimsup_{n \to \infty} A_n)=0$, where $\displaystyle \varlimsup_{n \to \infty} A_n:= \bigcap_{N=1}^\infty
\bigcup_{n=N}^\infty A_n$. Therefore
$$
\lim_{N \to \infty} \Pnum (\bigcup_{n=N}^\infty A_n)=0.
$$
Define $n_k$ inductively as below: $n_0:=0$ and
$$
n_k:=\inf\{N >n_{k-1}: \Pnum (\bigcup_{n=N}^\infty A_n) \leq \frac{1}{2^{k}}\}, \quad k=1,2, \cdots.
$$
Put $\displaystyle B_k:=\bigcup_{j \in (n_{k-1}, n_k]} A_j$. Then clearly
$$
\Pnum (B_k) \leq \Pnum (\bigcup_{j=n_{_{k-1}}+1}^{\infty} A_j) \leq \frac{1}{2^{k-1}},
$$
which implies \eqref{eq: prob-sum-B-finite}. \qed

\noindent \textbf{Proof of Theorem \ref{thm: main-B}.\;} Suppose there exists a subsequence $\{n_k\}_{k=1}^\infty$ of natural numbers such that \eqref{eq: prob-sum-B-infinite} and \eqref{eq: limit4Var-Expect}
hold true with $\displaystyle B_k:=\bigcup_{j \in (n_{k-1}, n_k]} A_j$ and $S_m:=1_{B_1}+\cdots+1_{B_m}$. Then Kochen Stone theorem tells us $\displaystyle \sum_{k=1}^\infty 1_{B_k} =\infty$ almost surely. And
\eqref{eq: equvi-infinite} tells us the validity of \eqref{eq: a.s.-infinite}.

Now assume the validity of \eqref{eq: a.s.-infinite}. It is just $\displaystyle \Pnum (\varlimsup_{n \to \infty} A_n)=1$, i.e.,
$$
\lim_{N \to \infty} \Pnum (\bigcup_{n=N}^\infty A_n)=\lim_{N \to \infty} \lim_{M \to \infty}\Pnum (\bigcup_{n=N}^M A_n)=1.
$$
We define $n_k$ inductively as the following: $n_0:=0$ and
$$
n_{k}:=\inf\{M > n_{k-1}: \Pnum (\bigcup_{n=n_{k-1}+1}^M A_n)>1-\frac{1}{2^{k}}\}, \quad k=1, 2, \cdots.
$$
Now put $\displaystyle B_k:=\bigcup_{j \in (n_{k-1}, n_k]} A_j$ and $S_m:=1_{B_1}+\cdots+1_{B_m}, \forall m \geq 1$. Clearly we have
$$
\Pnum (B_k)>1-\frac{1}{2^k}, k=1,2, \cdots.
$$
So $\Enum S_m >m -1$.

On the other hand
\begin{eqnarray*}
\Var (S_m) &=& \sum_{1 \leq k, \ell \leq m} \Cov (1_{B_k}, 1_{B_\ell}) \\
&\leq& \sum_{1 \leq k, \ell \leq m} \sqrt{\Var(1_{B_k}) \Var(1_{B_\ell})} \\
&\leq&  \sum_{1 \leq k, \ell \leq m} (\frac{1}{2})^{\frac{k+\ell}{2}} \leq \Bigl[\frac{\sqrt{1/2}}{1-\sqrt{1/2}}\Bigr]^2=3+2 \sqrt{2}.
\end{eqnarray*}
Now \eqref{eq: prob-sum-B-infinite} and \eqref{eq: limit4Var-Expect} are obvious.

\section{Discussions}\label{sec: C3}
In an old version of the current paper, the first author conjectured the following: Eq. \eqref{eq: a.s.-infinite} holds true if and only if there exists an increasing subsequence $\{n_k\}_{k=1}^\infty$ of
natural numbers such that \eqref{eq: limit4Var-Expect} and \eqref{eq: prob-sum-B-infinite} hold true with $B_k:=A_{n_k}$ and $S_m :=1_{B_1} +\cdots +1_{B_m}$. Afterwards, the second author present a
counter-example. And we reformulate it as below.

Let $X \sim U(\{1,2,3\})$, i.e. $X$ is valued in $\{1,2,3\}$ with equi-probability. Let
$$
A_{2n-1} :=\{X \neq 1\}, \quad A_{2n} :=\{X \neq 2\}
$$
for all natural number $n$. Clearly this sequence of events $\{A_n\}_{n=1}^\infty$ satisfies \eqref{eq: a.s.-infinite}. Then for any increasing subsequence $\{n_k\}_{k=1}
^\infty$ of natural numbers, define $B_k :=A_{n_k}$ and $S_m:=1_{B_1}+\cdots +1_{B_m}$, one clearly has $\Enum S_m=\frac{2m}{3}$. We claim further that,
\begin{equation}\label{eq: non-zero-limit}
\varliminf_{n \to \infty} \frac{\Var (S_n)}{(\Enum S_n)^2} \geq \frac{1}{8}>0.
\end{equation}
In fact, let
$$
t_m:=\sum_{k=1}^m 1_{\{n_k \hbox{ is odd}\}},
$$
then we can rewrite
\begin{eqnarray*}
S_m &=& t_m \cdot 1_{\{X \neq 1\}} +(m-t_m) \cdot 1_{\{X \neq 2\}}\\
&=& m-t_m \cdot 1_{\{X = 1\}} -(m-t_m) \cdot  1_{\{X = 2\}}.
\end{eqnarray*}
So we have
\begin{eqnarray*}
\Var (S_m) &=& \Var (t_m \cdot 1_{\{X = 1\}} +(m-t_m) \cdot  1_{\{X = 2\}})\\
&=& \frac{2}{9} \cdot [t_m^2+(m-t_m)^2 -t_m (m-t_m)] \geq \frac{m^2}{18}.
\end{eqnarray*}
Hence \eqref{eq: non-zero-limit} holds true.

\noindent \textbf{Acknowledgement.\;} The question discussed in this paper arises all in a sudden when the first author Xie was lecturing a probability course this semester. And he would like to thank his students for
their tolerance of his wandering to such questions during lessons. The observations \eqref{eq: equvi-finite} and \eqref{eq: equvi-infinite} are in fact inspired by Prof. Shihong Cheng (1939--2007) of Peking University
who taught Xie the course `Advanced Probability' in the fall of 1998. He appreciate Prof. Cheng's rigorism and carefulness in the course very much. This paper also serves as Xie's memoire of Prof. Cheng.


\begin{thebibliography}{00}
\bibitem{Chandra2012} Chandra, T. K.:
{\it The Borel-Cantelli Lemmas}, Springer Briefs in Statistics,
Springer, Heidelberg, 2012, xii+106 pp. ISBN: 978-81-322-0676-7

\bibitem{CK2001} Chernov, N.; Kleinbock, D.:
{\em Dynamical Borel–Cantelli lemmas for Gibbs measures}. 
Israel J. Math. \textbf{122} (2001), 1–-27.

\bibitem{CE1952} Chung, K. L.; Erdös, P.:
{\em On the application of the Borel-Cantelli lemma}.
Trans. Amer. Math. Soc. \textbf{72} (1952), 179–-186.

\bibitem{Davie2001} Davie, G.:
{\em The Borel-Cantelli lemmas, probability laws and Kolmogorov complexity}.
Ann. Probab. \textbf{29} (2001), no.4, 1426–-1434.

\bibitem{DMR2022} Dedecker, J.; Merlev\`{e}de, F.; Rio, E.:
{\em Criteria for Borel-Cantelli Lemmas with Applications to Markov Chains and Dynamical Systems},
Chapter 7 of {\it S\'{e}minaire de Probabilit\'{e}s LI}, Lecture Notes in Mathematics \textbf{2301}.

\bibitem{DF1965} Dubins, Lester E.; Freedman, David A.:
{\em A sharper form of the Borel-Cantelli lemma and the strong law}.
Ann. Math. Statist. \textbf{36} (1965), 800–-807.

\bibitem{Durrett} Durrett, R.:
{\it Probability: Theory and Examples},
Fifth edition. Cambridge Series in Statistical and Probabilistic Mathematics.
Cambridge University Press, Cambridge, 2019. xii+419 pp.

\bibitem{ER} Erd\"{o}s, P.; Rennyi, A.:
{\em On Cantor's series with convergent $1/q_n$},
Ann. Univ. Sci. Budapest E\"{o}tv\"{o}s Sect. Math. \textbf{2} (1959), 93--109.

\bibitem{EH2022} Estrada, L. F.; H\"{o}gele, M. A.:
{\em Moment estimates in the first Borel–Cantelli Lemma with applications to mean deviation frequencies},
Stat. Probab. Lett. \textbf{190} (2022), Paper No. 109636.

\bibitem{Freedman1973} Freedman, David:
{\em Another note on the Borel-Cantelli lemma and the strong law, with the Poisson approximation as a by-product}.
Ann. Probability \textbf{1} (1973), 910–-925.

\bibitem{Frolov2021} Frolov, A. N.:
{\em On strong forms of the Borel–Cantelli lemma and intermittent interval maps},
J. Math. Anal. Appl. \textbf{504} (2021), Paper No. 125425.

\bibitem{Hoover2002} Hoover, Donald R.:
{\em First Borel-Cantelli equalities incorporating known dependency structure}.(English, Italian summary)
Metron \textbf{60} (2002), no.3-4, 21-–30.

\bibitem{HLSW2022} Hussain, M.; Li, B.; Simmons, D.; Wang, B.:
{\em Dynamical Borel–Cantelli lemma for recurrence theory},
Ergod. Th. \& Dynam. Sys. \textbf{42} (2022), 1994–-2008.

\bibitem{Kifer2021} Kifer, Y.:
{\em The Strong Borel–Cantelli Property in Conventional and Nonconventional Setups},
Lecture Notes in Math.Vol \textbf{2290}. CIRM Jean-Morlet Ser. Springer, Cham, 2021, 235–-261.
ISBN: 978-3-030-74862-3; 978-3-030-74863-0


\bibitem{KS1964} Kochen, S. B., Stone, C. J.:
{\em A note on the Borel-Cantelli lemma}. 
Illinois J. Math. \textbf{8} (1964), 248--251. 

\bibitem{Lamperti1963} Lamperti, J.: 
{\em Wieners test and Markov chains}. 
J. Math. Anal. Appl. \textbf{6} (1963), 58--66.

\bibitem{Luzia2014} Luzia, N.:
{\em A Borel-Cantelli lemma and its applications},
Trans. Amer. Math. Soc. \textbf{366} (2014), no.1, 547-–560.

\bibitem{Petrov2002} Petrov, Valentin V.:
{\em A note on the Borel-Cantelli lemma}.
Statist. Probab. Lett. \textbf{58} (2002), no.3, 283–-286.

\bibitem{Philipp1967} Philipp, W.:
{\em Some metrical theorems in number theory}. 
Pacific J. Math. \textbf{20} (1967), 109–-127.

\bibitem{Simon} Simon, B.:
{\it Functional Integration and Quantum Physics},
2nd Edition. AMS Chelsea Publishing, American Mathematical Society, Providence,
Rhode Island.

\bibitem{XZ} Xie, Jiansheng; Zhou, Ruisong:
{\em A note on an SLLN of Dvoretzky-Erdös type}.
Statist. Probab. Lett. \textbf{195} (2023), Paper No. 109769, 5 pp.


\end{thebibliography}
\end{document}